\documentclass{cimart}

\usepackage{bbm}
\usepackage[all]{xy}

\DeclareMathOperator{\ob}{Ob}

\newcommand{\aut}{\mathrm{Aut}}
\newcommand{\calc}{\mathcal{C}}
\newcommand{\dg}[1]{^{[#1]}}
\newcommand{\f}{\mathcal{F}}
\newcommand{\fb}{\mathbf{FB}}
\newcommand{\fin}{\mathbf{FA}}
\newcommand{\finj}{\mathbf{FI}}
\newcommand{\fs}{\mathbf{FS}}
\newcommand{\kmod}{\mathtt{Mod}_\kring}
\newcommand{\kring}{\mathbbm{k}}
\newcommand{\nat}{\mathbb{N}}
\newcommand{\op}{^\mathrm{op}}
\newcommand{\pbar}{\overline{P}}
\newcommand{\primfs}{{\kring \fs^0}}
\newcommand{\sgn}{\mathsf{sgn}}
\newcommand{\sym}{\mathfrak{S}}
\newcommand{\triv}{\mathsf{triv}}

\theoremstyle{cmdefinition}
\newtheorem{notation}[definition]{Notation}

\title{Filtering the linearization of the category of surjections}

\author{Geoffrey Powell}

\authorinfo{University of Angers, France}{Geoffrey.Powell@math.cnrs.fr}

\abstract{%
    A filtration of the morphisms of the $\kring$-linearization $\kring \fs$ of the category $\fs$ of finite sets and surjections is constructed using a natural $\kring \finj\op$-module structure induced by restriction, where $\finj$ is the category of finite sets and injections. In particular, this yields the `primitive' subcategory $\primfs \subset \kring \fs$ that is of independent interest; for example, the category of $\primfs$-modules is closely related to the category of $\kring \fin$-modules, where $\fin$ is the category of finite sets and all maps. Working over a field of characteristic zero, the subquotients of this filtration are identified as bimodules over $\kring \fb$, where $\fb$ is the category of finite sets and bijections, also exhibiting and exploiting additional structure. In particular, this describes the underlying $\kring \fb$-bimodule of $\primfs$.
    }



\VOLUME{34}
\YEAR{2026}
\ISSUE{1}
\NUMBER{5}
\DOI{https://doi.org/10.46298/cm.16066}

\begin{document}

\section{Introduction}
\label{sect:intro}

The main player of this paper is  $\kring \fs$, the $\kring$-linearization of the category $\fs$ of finite sets and surjections, where $\kring$ is a field of characteristic zero. Representations of this category (known as $\kring\fs$-modules here) and of its opposite are of significant interest. For instance, $\kring \fs \op$-modules play an important role in the application of Sam and Snowden's Gröbner-theoretic methods \cite{MR3556290};  Tosteson \cite{MR4546495} has provided additional interesting information on $\kring \fs\op$-modules. 

In addition, Pirashvili \cite{Phh} proved that the category of $\kring \fs$-modules is Morita equivalent to that of $\kring \Gamma$-modules, where $\Gamma$ is the category of finite pointed sets and pointed maps (this does not require working over a field of characteristic zero). The category of $\kring \Gamma$-modules occurs in numerous situations; for instance, Pirashvili exhibited $\kring \Gamma$-modules as the natural coefficients for his higher Hochschild homology when working with pointed spaces. For {\em unpointed} spaces, $\kring \Gamma$ should be replaced by $\kring \fin$, where $\fin$ is the category of finite sets and all maps. The interplay between these structures was investigated in the joint work with Vespa \cite{MR4518761}, motivated by the joint work \cite{PV} on higher Hochschild homology evaluated on wedges of circles. The latter was inspired by the work of Turchin and Willwacher \cite{MR3982870}.

The purpose of the current paper is  to construct a (highly non-trivial) filtration of the morphism spaces of $\kring \fs$ and to investigate its properties. 
 This is based  upon the embedding $\fs \subset \fin$ and the following observation: for a finite set $X$, there is a natural $\kring \fin\op$-module structure on $
\kring \fs (- , X)$ such that the projection $\kring \fin (-, X) \rightarrow \kring \fs (-, X)$ that sends non-surjective maps to zero is a morphism of $\kring \fin\op$-modules.

It follows that $\kring \fs (- , X)$ has a $\kring \finj\op$-module structure given by restriction along $\kring \finj \op\hookrightarrow \kring \fin\op$, where $\finj \subset \fin$ is the category  of finite sets and injections.  General properties of $\kring \finj\op$-modules then yield  the exhaustive filtration: 
\begin{eqnarray}
\label{eqn:filt_fs_intro}
0=\kring \fs ^{-1}(-, X) \subseteq \primfs (-, X) \subseteq \kring \fs^1 (-, X)  \subseteq\kring \fs^2  (-, X ) \subseteq\ldots \subseteq \kring \fs (-, X)
\end{eqnarray}
(see Corollary \ref{cor:filt_qfs}). 

As an example, for a finite set $U$,  $\primfs (U, X)$ is the subspace of $\kring \fs (U,X)$ given by the intersection of the kernels of the $\kring$-linear  maps 
$$
\kring \fs (U,X)
\rightarrow \kring \fs (V, X)
$$
given by the $\kring \finj\op$-module structure of $\kring \fs (-, X)$,  for non-bijective injections $V \hookrightarrow U$. Whilst elementary to define, determining $\primfs (U,X)$ as a $\kring$-vector space or, better still, a $\kring ( \aut (U) \op\times\aut (X) )$-module  (where $\aut (-)$ denotes the  automorphism group), is not easy.

There is more structure,  which shows the significance of $\primfs$:

\begin{theorem}
\label{THM:primfs}
(Corollary \ref{cor:properties_filtration}.) The following properties hold:
\begin{enumerate}
\item 
$\primfs$ forms  a wide $\kring$-linear subcategory of $\kring \fs$; 
\item 
$\kring \fs$ has the structure of a $\primfs \otimes \kring \finj\op$-module and  (\ref{eqn:filt_fs_intro}) is a filtration by $\primfs \otimes \kring \finj\op$-modules.
\end{enumerate}
\end{theorem}

Part of the interest of the $\kring$-linear category $\primfs$ is that the category of $\primfs$-modules is closely related to that of $\kring \fin$-modules. The study of this relationship was initiated in \cite{MR4518761} and developed in \cite{P_finmod}: a precise relationship with the category of $\kring \fin$-modules is given by the Morita equivalence result \cite[Theorem 7.10]{P_finmod}, which is expressed using a $\kring$-linear category that is closely related to $\primfs$. 

\begin{remark}
A  conceptual explanation for the relationship between $\primfs$ and $\kring \fin$ is provided in \cite[Section 19]{P_rel_kos} using relative nonhomogeneous Koszul duality {\em à la Positselski} \cite{MR4398644}. This is based on the fact that every map between finite sets can be written as the composite of a surjection followed by an injection; this  induces an identification of the $\kring$-linear category  $\kring \fin$ as $\kring \finj \otimes_{\fb} \kring \fs$, equipped with composition defined using a distributive law. Relative non-homogeneous Koszul duality then yields a DG $\kring$-linear category which has underlying complex that identifies with  the Koszul complex that is used in \cite{MR4518761}. (This approach is not exploited explicitly here.)
\end{remark}

To study the filtration (\ref{eqn:filt_fs_intro}), we consider the filtration subquotients $\kring \fs^{\ell /(\ell-1)}$, for each $\ell \in \nat$ (for $\ell =0$, this subquotient is simply $\primfs$). By Theorem \ref{THM:primfs}, these subquotients are $\primfs \otimes \kring \finj\op$-modules and the aim is to understand this structure. 
The first step is to give a precise relationship between the higher filtration quotients and $\primfs$.

In Theorem \ref{THM:subquotients} below, $\sgn_a$ denotes the sign representation of the symmetric group $\sym_a$ and $S_{(\ell, 1^a)}$ the simple representation of $\sym_{\ell +a}$ indexed by the partition $(\ell , 1^a)$ (considered as a $\kring \sym_{\ell +a} \op$-module). Then $\bigoplus_{a\in \nat}  \sgn_a \boxtimes S_{(\ell, 1^a)} $ is considered as a $\kring \fb \otimes \kring \fb\op$-module. Moreover, $\odot_{\fb\op}$ denotes the Day convolution product for $\kring \fb\op$-modules (see Definition~\ref{defn:day_convolution})  and $\triv_\ell$ is the trivial $\kring \sym_\ell\op$-representation, considered as a $\kring \fb\op$-module. 

We also use that $\kring \fs$  is equipped with an augmentation $\kring \fs \rightarrow \kring \fb$ that sends non-bijections to zero. Tensoring $\primfs \rightarrow \kring \fb$,  the restriction of the above augmentation to $\primfs$,   with $\kring \fb\op$ yields the $\kring$-linear functor  $\primfs \otimes \kring \fb\op \rightarrow \kring \fb\otimes \kring \fb\op $.

\begin{theorem}[{Theorem \ref{thm:coker_theta_ell}}]
\label{THM:subquotients}
For $\ell >0$, there is a short exact sequence of $\primfs \otimes \kring \fb\op$-modules
\[
0
\rightarrow 
\kring \fs^{\ell/(\ell -1)}
\stackrel{\theta^\ell}{\longrightarrow} 
\primfs \odot_{\fb\op} \triv_\ell
\rightarrow 
\bigoplus_{a\in \nat}  \sgn_a \boxtimes S_{(\ell, 1^a)} 
\rightarrow 
0,
\]
where the cokernel is a $\primfs \otimes \kring\fb\op$-module via  restriction along the homomorphism $\primfs \otimes \kring \fb\op \rightarrow \kring \fb\otimes \kring \fb\op $.
\end{theorem}

The way in which the cokernel of  $\theta^\ell$ appears is of interest; to explain this,  we use the skeleton of the category $\fin$ (and hence its wide subcategories) provided by the sets $\mathbf{n} = \{1, \ldots, n \}$, for $n \in \nat$. Taking $b = \ell + a $ and evaluating $\theta^\ell$ on $(\mathbf{b}, \mathbf{a})$, the map identifies as a morphism of $\kring (\sym_a \times \sym_b\op)$-modules:
\[
\Theta _\mathbf{a} (\mathbf{b}) : \kring \fs (\mathbf{b}, \mathbf{a}) \rightarrow  \kring ^{\finj (\mathbf{a}, \mathbf{b})}. 
\]
 The cokernel of this map  identifies as $\sgn_a \boxtimes S_{(\ell, 1^a)} $, by Theorem \ref{thm:coker_Theta} and Lemma \ref{lem:identify_Lambda^a_pbar}. 
 Then, by exploiting the $\primfs$-module structure, it is shown that, as $a$ varies, this accounts for all of the cokernel of $\theta^\ell$.

Theorem \ref{THM:subquotients} gives an explicit relationship between $\kring \fs$ and $\primfs \odot_{\fb\op} \triv$, where $\triv$ is the constant $\kring \fb\op$-module $\kring$, which encodes all of the $\triv_\ell$. Working up to isomorphism in the category of  $\kring \fb\op$-modules, $-\odot_{\fb\op} \triv$ is invertible, which leads to:

\begin{theorem}[{Theorem~\ref{thm:primfs_bimodule}}] \label{THM:primfs_identify}
In the Grothendieck group of $\kring \fb \otimes \kring \fb\op$-modules, there is an equality
\begin{eqnarray*}
[\primfs] + \sum_{b> a \geq 0} (-1)^{b-a} [\sgn_a \boxtimes \sgn_b] 
= 
[\kring \fs] \odot_{\fb\op} \Big( 
\sum_{t \in \nat} (-1)^t  [\sgn_t]
\Big) .
\end{eqnarray*}
\end{theorem}

This result is equivalent to \cite[Theorem 3]{MR4518761} which was established using  different methods (see Remark \ref{rem:compare_Theorem3_PV2} for a brief explanation of the relationship). The current approach is designed to give a quick derivation of this, working over a field of characteristic zero, that also determines the higher filtration quotients. 

As already indicated, the work of  \cite{MR4518761} was motivated in part by the study of the relationship between $\primfs$-modules and the category of $\kring \fin$-modules.  The subtlety   arises from the exceptional behaviour of the simple $\kring\fin$-modules  indexed by the partitions of the form $(1^t)$. This is behind the appearance of the terms $\sgn_a$ in the statements of Theorems~\ref{THM:subquotients} and~\ref{THM:primfs_identify}.

This exceptional behaviour is already apparent in the classification of the simple $\kring \fin$-modules given by Wiltshire-Gordon in \cite{2014arXiv1406.0786W}.  (This classification is outlined very briefly in Section \ref{subsect:coker_Theta} and is revisited in \cite{P_finmod}.) Moreover, Sam, Snowden and Tosteson \cite{MR4788713}  observed that the classification is related to that of certain simple representations of the Witt Lie algebra by Feigin and Shoikhet \cite{MR2350124}, where this exceptional behaviour is also apparent (see Remark \ref{rem:Feigin-Shoikhet} for slightly more detail).

\begin{remark}
Throughout this paper, $\kring$ is taken to be a field of characteristic zero. Given that many of the results of \cite{MR4518761} hold over an arbitrary commutative ring, it is natural to seek to generalize the current results.  Indeed, the construction of the filtration (\ref{eqn:filt_fs_intro}) can be carried out over any such ring and Theorem \ref{THM:primfs} should also hold in this generality.  

However, the proof of Theorem \ref{THM:subquotients} depends upon $\kring$ being a field of characteristic zero, since it is based upon Theorem \ref{thm:coker_Theta}, the given proof of which exploits the behaviour of the simple $\kring \fin$-modules over a field of characteristic zero. It would  be of considerable interest to have a characteristic-free version of Theorem \ref{thm:coker_Theta}.
\end{remark}

\section{Background}
\label{sect:background}

This section  presents some background and introduces the main characters, namely some familiar wide subcategories of the category of finite sets. The main objective is to state and sketch the proof of Theorem \ref{thm:coker_Theta} (which is proved in detail in \cite[Section 6]{P_finmod}); this  is the initial input for the proofs of the main theoretical results of the paper in Section \ref{sect:subquotients}.

Throughout  $\kring$ is a field  of characteristic zero. Tensor products of $\kring$-vector spaces are denoted simply by $\otimes$.

\subsection{Basics}
\label{subsect:basics}

The symmetric group on $b$ letters is denoted by  $\sym_b$. The hypothesis upon $\kring $  allows us to invoke the (semi-simple) representation theory of the symmetric groups in characteristic zero. (Basics of the representation theory of the symmetric groups in characteristic zero are covered, for example, in \cite[Chapter 4]{MR1153249}.) 
 The simple representations of $\sym_b$  are indexed by the partitions $\lambda \vdash b$, with  the corresponding simple module denoted by $S_\lambda$. We use the convention that the partition $(b)$ indexes the trivial representation,  $ \triv_b= S_{(b)}$.

For a category $\calc$, the morphism sets are denoted either $\hom_\calc (-,-)$ or $\calc (-,-)$. 
The $\kring$-linearization of $\calc$ is denoted by $\kring \calc$; this has the same objects as $\calc$ and morphisms given by $ \hom_{\kring \calc} (-,-) := \kring \hom_\calc (-,-)$, also written $\kring \calc (-,-)$. For $f$ a morphism of $\calc$, the corresponding morphism of $\kring \calc$ is denoted $[f]$; such elements provide  $\kring$-bases for the morphism spaces. Moreover, this defines a faithful functor $\mathcal{C} \hookrightarrow \kring \mathcal{C}$ that is the identity on objects.

If $\calc$ is essentially small, one can consider the category of functors from $\calc$ to $\kmod$, the category of  $\kring$-vector spaces; 
this is sometimes denoted $\f (\calc)$.  This is equivalent to the category of $\kring$-linear functors from $\kring \calc$ to $\kring$-vector spaces, which will  be referred to here as the category of $\kring \calc$-modules (rather than $\calc$-modules, as is common in the literature, since we also consider $\kring$-linear categories that are not of the form $\kring \calc$). 

The category $\f (\calc)$ is abelian and comes equipped with the tensor product induced by that on $\kring$-vector spaces; namely, for functors $F, G$, their tensor product is defined by $(F \otimes G)(X):= F(X) \otimes G(X)$. Vector space duality induces  $D : \f (\calc ) \op \rightarrow \f (\calc \op)$ given by $DF (X) := \hom_\kring (F(X), \kring)$. This restricts to an equivalence between the respective full subcategories of functors taking finite-dimensional values. 

The tensor product $\otimes _\calc : \f (\calc\op) \times \f(\calc) \rightarrow \kmod$ is defined as the composite of the external tensor product $\boxtimes : \f(\calc\op) \times \f(\calc) \rightarrow \f (\calc \op \times \calc)$ (which is simply given by $(F \boxtimes G) (X \times Y) = F(X) \otimes G(Y))$ and the coend $\f (\calc \op \times \calc) \rightarrow \kmod$.

A subcategory $\calc'$ of $\calc$ is {\em wide} if it contains all the objects of $\calc$; a wide subcategory is  determined by specifying the morphisms that it contains.

\subsection{The cast}
\label{subsect:cast}

The category $\fin$ of finite sets and {\em all} maps has skeleton $\{ \mathbf{n} \ | \ n \in \nat\}$, where $\mathbf{n}$ is the set $\{1, \ldots, n \}$ (understood to be $\emptyset$ for $n=0$). An obvious but useful fact is that, for finite sets $X$ and $Y$, $\fin(X, Y)$ is a finite set.

The  wide subcategories $\fs$, $\finj$, $\fb$ of $\fin$ are our main players, specified by their  classes of morphisms:  
 surjections, injections, bijections respectively. These fit into the commutative diagram of inclusions of wide subcategories:
\[
\xymatrix{
\fb
\ar@{^(->}[r]
\ar@{^(->}[d]
&
\finj
\ar@{^(->}[d]
\\
\fs
\ar@{^(->}[r]
&
\fin.
}
\]
Any map of $\fin$ factors (uniquely up to a morphism of $\fb$) as the composite of a map of $\fs$ followed by one of $\finj$.

\begin{remark}
\label{rem:fb_groupoid}
Since $\fb$ is a groupoid, it is isomorphic to its opposite, $\fb\op$. Thus the category of $\kring \fb$-modules is isomorphic to the category of $\kring \fb\op$-modules. This is a global form of the fact that, for $b \in \nat$, the category of $\kring \sym_b$-modules is isomorphic to the category of $\kring \sym_b\op$-modules. This will be used to adjust variance, when necessary, without further comment.
\end{remark}

If $f: X \rightarrow Y$ is a map between finite sets that is not surjective, for any map $g : U \rightarrow X$, precomposing with $g$ gives $f \circ g$ that is not surjective. Likewise, if $f$ is not injective, then postcomposing with any map $h: Y \rightarrow Z$  gives $h \circ f$ that is not injective. This implies the following:

\begin{proposition}
\label{prop:quotient_fs_finj}
Let $X$ be a finite set. 
\begin{enumerate}
\item 
There is a unique $\kring \fin$-module structure on $\kring \finj (X, -)$ extending the canonical $\kring \finj$-module structure such that the projection 
\[
\kring \fin (X, - ) \twoheadrightarrow \kring \finj (X,-)
\]
sending non-injective maps to zero is a morphism of $\kring \fin$-modules.  
\item 
There is a unique $\kring \fin\op$-module structure on $\kring \fs (-, X)$ extending the canonical $\kring \fs\op$-module structure such that the projection 
\[
\kring \fin (-, X ) \twoheadrightarrow \kring \fs (-,X)
\]
sending non-surjective maps to zero is a morphism of $\kring \fin\op$-modules.  
\item 
Evaluating on $X$, these both yield the surjective morphism of $\kring$-algebras 
\[
\kring \fin (X, X) 
\twoheadrightarrow 
\kring \aut (X)
\]
induced by sending non-bijective maps to zero. This is split by the canonical inclusion of $\kring$-algebras $\kring \aut (X) \hookrightarrow \kring \fin (X,X)$. 
\end{enumerate}
\end{proposition}

For finite sets $X, U ,Y$, the surjection $\kring \fin (-, X) \twoheadrightarrow \kring \fs (-,X)$ yields the following commutative diagram corresponding to the $\kring \fin\op$-module structures: 
\[
\xymatrix{
\kring \fin (U, X) \otimes \kring \fin (Y,U) 
\ar[r]
\ar@{->>}[d]
&
\kring \fin (Y, X) 
\ar@{->>}[d]
\\
\kring \fs (U, X) \otimes \kring \fin (Y,U) 
\ar[r]
&
\kring \fs (Y, X). 
}
\]
Taking $X=Y$, this refines to 
\[
\xymatrix{
\kring \fin (U, X) \otimes \kring \fin (X,U) 
\ar[r]
\ar@{->>}[d]
&
\kring \fin (X, X) 
\ar@{->>}[d]
\\
\kring \fs (U, X) \otimes \kring \finj (X,U) 
\ar[r]
&
\kring \aut (X), 
}
\]
using  the projection $\kring \fin (X, U) \twoheadrightarrow \kring \finj (X, U)$ to yield the left hand vertical arrow.

The lower horizontal map is compatible with the $\kring \fin$-actions on $\kring \fs (-,X)$ and $\kring \finj (X, -)$ respectively, yielding the following: 

\begin{proposition}
\label{prop:fs_fin_finj_composition}
For a finite set $X$, composition in $\fin$ induces the morphism of $\kring \aut (X)$-bimodules:
\begin{eqnarray}
\label{eqn:fs_finj}
\kring \fs (-, X) \otimes_\fin \kring \finj (X, -) 
\rightarrow 
\kring \aut (X). 
\end{eqnarray}
\end{proposition}

\begin{proof}
Associativity of composition implies that, for a finite set $U$, the morphism 
\[
\kring \fs (U, X) \otimes \kring \finj (X, U) 
\rightarrow 
\kring \aut (X)
\]
induced by composition and projection to $\kring \aut (X)$ is a morphism of $\kring \aut (X)$-bimodules for the obvious bimodule structures.

It  remains to check that these induce a morphism $\kring \fs (-, X) \otimes_\fin \kring \finj (X, -) 
\rightarrow 
\kring \aut (X)$ as claimed. This follows from the associativity of the composition in $\kring \fin$ together with the fact that the factorization observed above across the surjection of $\kring \fin$-modules $\kring \fin (X, -) \twoheadrightarrow \kring \finj (X, -)$ is canonical.

It is worth spelling this out. For composable morphisms $f \in \fs (U, X)$, $g \in \fin (V,U)$ and $h \in \finj (X, V)$, we require to prove that the images of $g^* [f] \otimes [h]$ and $[f] \otimes g_* [h]$ in $\kring \aut (X)$ are the same, writing $g^*$ and $g_*$ to indicate the $\kring \fin\op$ and $\kring \fin$-actions respectively. Now, $g^*[f] = [f \circ g]$ if $f\circ g \in \fin (V, X)$ is surjective and zero otherwise; likewise $g_* [h]=[g \circ h]$ if $g \circ h \in \fin (X, U)$ is injective and zero otherwise. Similarly $ f\circ g \circ h \in \fin (X, X)$ has image $[f \circ g \circ h]$  in $ \kring \aut (X)$ if $f \circ g \circ h$ is a bijection and zero otherwise. Now, the  condition that $f \circ g \circ h$ is a bijection implies that  $f \circ g$ is surjective and that $g \circ h$ is injective. 
 It follows  that the images of  $g^* [f] \otimes [h]$ and $[f] \otimes g_* [h]$ in $\kring \aut (X)$ are  both equal to $[f \circ g \circ h]$ if $f \circ g \circ h$ is a bijection and both zero otherwise.
\end{proof}

We use this Proposition to construct the natural transformation $\Theta _X$ of equation (\ref{eqn:map_fs_finj}). This is used in Section \ref{sect:subquotients},  where it is essential input in the proof of Theorem \ref{thm:coker_theta_ell}, via Theorem \ref{thm:coker_Theta}, the main result of this section,  which identifies its cokernel.

\begin{corollary}
\label{cor:adjoint_X}
For a finite set $X$, 
\begin{enumerate}
\item 
by adjunction, (\ref{eqn:fs_finj})  yields the morphism of $\kring \fin \op \otimes \kring \aut (X)$-modules:
\begin{eqnarray}
\label{eqn:fs_adjoint_X=Y}
\kring \fs (-, X) \rightarrow \hom_{\kring \aut (X)\op} (\kring \finj (X, -), \kring \aut (X));
\end{eqnarray}
\item 
the codomain of this morphism is isomorphic to $D \kring \finj (X, -)$ equipped with its canonical $\kring \fin \op \otimes \kring \aut (X)$-module structure, hence (\ref{eqn:fs_adjoint_X=Y}) identifies as the morphism 
\begin{eqnarray}
\label{eqn:map_fs_finj}
\Theta_X : \kring \fs (-, X) \rightarrow D \kring \finj (X, - ) = \kring ^{\finj (X,-)}
\end{eqnarray}
of $\kring \fin \op \otimes \kring \aut (X)$-modules.
\item 
For a finite set $U$, $\Theta_X (U)$ identifies as 
\begin{eqnarray*}
\kring \fs (U, X) 
&\rightarrow 
&
\kring ^{\finj (X,U)}
\\
\ [f]& \mapsto & \sum_{s\in \mathrm{Sec}(f)} \eta_s,
\end{eqnarray*}
 where $\mathrm{Sec}(f) \subset \finj (X, U)$ is the set of sections of $f \in \fs (U,X)$ and $\eta_s$ is the dual basis element to $[s] \in \kring \finj (X,U)$.
\end{enumerate}
\end{corollary}

\begin{proof}
The first statement is given by the standard adjunction, taking into account the tensor product $\otimes _\fin$ and the $\kring \aut (X)$-bimodule structure. 

The second statement follows from the general fact that, for a finite group $G$ and a right $\kring G$-module $M$, postcomposing with the $\kring$-linear map $\kring G \rightarrow \kring $ that sends $[e] \mapsto 1$ and $[g] \mapsto 0$ for $g \neq e$ yields an isomorphism of left $\kring G$-modules:
\[
\hom_{\kring G\op} (M, \kring G) 
\rightarrow 
\hom_\kring (M, \kring) ,
\]
where $\kring G$ is considered as a $\kring G$-bimodule with respect to the left and right regular actions. 
 The inverse sends a linear map $\phi \in \hom_\kring (M, \kring)$ to the linear map $M \rightarrow \kring G$ defined by $x \mapsto \sum_{g \in G} \phi (x g) \otimes [g^{-1}]$.  Since this is natural with respect to $M$, the second assertion follows.

The explicit identification of $\Theta_X$ follows immediately from this.
\end{proof}

\begin{example}
\label{exam:indentify_Theta}
In the case $U=X$, $\Theta_X (X)$ is the linear map $\kring \aut (X) \rightarrow \kring ^{\aut (X)}$ given by $[\alpha]\mapsto \eta_{\alpha^{-1}}$. This is  an isomorphism.
\end{example}

One can restrict the $\kring \fin\op$-module structure of $\kring \fs (-,X)$ to give a $\kring \finj\op$-module structure, with structure morphism of  $\kring \aut (X) \otimes \kring \aut (Y)\op$-modules: 
\begin{eqnarray*}
\kring \fs (U, X) \otimes \kring \finj (Y, U) 
\rightarrow 
\kring \fs (Y, X). 
\end{eqnarray*}
Moreover, these factors across  
\begin{eqnarray}
\label{eqn:XneqY_fs_finj}
\kring \fs (U, X) \otimes_{\aut (U)}  \kring \finj (Y, U) 
\rightarrow 
\kring \fs (Y, X).
\end{eqnarray}

The following is a counterpart of Corollary \ref{cor:adjoint_X}:

\begin{proposition}
\label{prop:adjoint_XY}
For finite sets $U, X,Y,$ (\ref{eqn:XneqY_fs_finj}) yields the morphism of $\kring \aut (U) \op \otimes \kring \aut (X)$-modules, by adjunction:
\[
\kring \fs (U, X)
\rightarrow 
\hom_{\kring \aut (Y) \op} 
(\kring \finj (Y, U) , \kring  \fs (Y, X) ).
\]

When $X=Y$, this identifies with  (\ref{eqn:fs_adjoint_X=Y}) evaluated on $U$, with the inherited $\kring \aut (U)\op$-action.
\end{proposition}

\subsection{The cokernel of $\Theta$}
\label{subsect:coker_Theta}

For compatibility with later notation, here we take $X$ to be $\mathbf{a}$ (for some $a \in \nat$) and consider the natural transformation (\ref{eqn:map_fs_finj}) 
\[
\Theta_\mathbf{a} : 
\kring \fs (-, \mathbf{a}) \rightarrow D \kring \finj (\mathbf{a}, -) 
\]
of $\kring \fin\op \otimes \kring \sym_a$-modules. 

The $\kring \fin\op$-module structure is useful here. Simple $\kring \fin$-modules are understood by the work of Wiltshire-Gordon \cite{2014arXiv1406.0786W}, in particular his classification  \cite[Theorem 5.5]{2014arXiv1406.0786W}. (This theory has been revisited in \cite[Section 6]{P_finmod}, with the classification stated as \cite[Theorem~6.6]{P_finmod}.) By dualizing, one obtains the classification of the simple $\kring \fin\op$-modules. We require a few details, as sketched below (for more details, see \cite[Section 6]{P_finmod}).

The standard projective $\kring \fin$-module $P_\mathbf{1}:= \kring \fin (\mathbf{1}, -)$  identifies with the functor $X \mapsto \kring X$. There is a natural transformation $P_\mathbf{1} \rightarrow \kring$ to the constant functor with value $\ kring $ induced by the canonical map $X \rightarrow *$ to the final object of $\fin$. The kernel is denoted $\overline{P}$, which lives in the short exact sequence 
\[
0
\rightarrow 
\overline{P}
\rightarrow 
P_\mathbf{1}
\rightarrow 
\overline{\kring}
\rightarrow 
0,
\]
where $\overline{\kring}\subset \kring$ is the subfunctor supported on non-empty sets. This short exact sequence does not split: both $\overline{P}$ and $\overline{\kring}$ are simple $\kring \fin$-modules, whereas  neither are projective.

For each $t \in \nat$, one can form the $\kring \fin$-module $\Lambda^t (\pbar)$ by post-composing with the $t$th exterior power functor; by convention, $\Lambda^0 (\pbar)$ is taken to be $\overline{\kring}$. For each $t$, $\Lambda^t (\pbar)$ is a  simple $\kring \fin$-module.  

The inclusion $\pbar \hookrightarrow P_\mathbf{1}$ induces $\Lambda^t (\pbar) \hookrightarrow \Lambda^t (P_\mathbf{1})$ which (for $t>0$) fits into the short exact sequence 
\[
0
\rightarrow 
\Lambda^t (\pbar) 
\hookrightarrow 
\Lambda^t (P_\mathbf{1}) 
\rightarrow 
\Lambda^{t-1} (\pbar)
\rightarrow 
0. 
\]

\begin{lemma}
\label{lem:identify_Lambda^a_pbar}
For $a \in \nat$, the underlying $\kring \sym_b$-module of $\Lambda^a (\pbar ) (\mathbf{b})$ is isomorphic to 
$S_{(b-a,1^a)}$ if $b>a$ and $0$ otherwise. 
\end{lemma}

\begin{proof}
The case $a=0$ is immediate, using the convention $\Lambda^0 (\pbar) = \overline{\kring}$. Hence, suppose that $a>0$.

First consider $\Lambda^a (P_\mathbf{1}) (\mathbf{b}) = \Lambda^a (\kring \mathbf{b})$. This is zero if $a>b$. Otherwise one checks that it identifies with the $\kring \sym_b$-module $(\sgn_b \boxtimes \triv_{b-a}) \uparrow_{\sym_b \times \sym_{b-a}} ^{\sym_b} \cong S_{(b-a,1^a)} \oplus S_{(b-a+1, 1^{a-1})}$, where the identification follows from the Pieri formula. (The Pieri formula is a special case of the Littlewood-Richardson rule; an equivalent version in terms of the  associated Schur functors is stated  in \cite[Section 1.1]{MR3430359}, for example.)

 The result then follows from an obvious induction on $a$ using the above short exact sequence.
\end{proof}

The following result is the dual of \cite[Proposition 6.11]{P_finmod}. The proof is sketched here, for completeness.

\begin{theorem}
\label{thm:coker_Theta}
For $a \in \nat$, the natural transformation $\Theta_\mathbf{a}$ fits into  an exact sequence of $\kring \fin\op \otimes \kring \sym_a$-modules: 
\[
\kring \fs (-, \mathbf{a})
\stackrel{\Theta_\mathbf{a}}{ \longrightarrow}
 D \kring \finj (\mathbf{a}, -) 
\rightarrow 
D \Lambda^a (\pbar) \boxtimes \sgn_a
\rightarrow 
0.
\]
In particular, the cokernel of $\Theta_\mathbf{a}$ is a simple $\kring \fin\op \otimes \kring \sym_a$-module.
\end{theorem}

\begin{proof}
The key point of the proof is the identification of $\kring \finj( \mathbf{a}, -)$ as a $\kring \fin \otimes \kring \sym_a\op$-module. 

Consider the surjection $\kring \fin (\mathbf{a}, -) \twoheadrightarrow \kring \finj (\mathbf{a}, -)$ of $\kring \fin$-modules given by Proposition~\ref{prop:quotient_fs_finj}. As a $\kring \fin$-module, $\kring \fin (\mathbf{a}, -)$ is isomorphic to $P_1^{\otimes a}$, hence the inclusion $\pbar \subset P_1$ induces $\pbar^{\otimes a} \subset \kring \fin (\mathbf{a}, -)$. Composing these gives the morphism of $\kring \fin \otimes \kring \sym_a\op$-modules:
\[
\pbar^{\otimes a} \rightarrow \kring \finj (\mathbf{a}, -).
\]

One shows (see \cite[Proposition 6.1]{P_finmod})  that the image of this map is semisimple and is isomorphic to 
\[
\bigoplus_{\substack{\lambda \vdash a \\ \lambda \neq (1^a) }} \big( \kring \finj (\mathbf{a}, -)\otimes_{\sym_a} S_\lambda\big) \boxtimes S_\lambda 
\quad 
\oplus 
\quad 
\Lambda^a (\pbar) \boxtimes \sgn_a.
\]
The cokernel of the map is the simple object  $
\Lambda^{a-1} (\pbar) \boxtimes \sgn_a$.

By \cite[Theorem 6.6]{P_finmod}, the simple $\kring \fin$-modules $\kring \finj (\mathbf{a}, -)\otimes_{\sym_a} S_\lambda$ ($\lambda \neq (1^a)$) and $\Lambda^{a-1} (\pbar)$ vanish on evaluation on $\mathbf{n}$ with $n <a$ and are detected by evaluation on $\mathbf{a}$. Likewise $\Lambda^a (\pbar)(\mathbf{n})=0$ for $n \leq a$ and $\Lambda^a (\pbar)(\mathbf{a+1})\neq 0$.  On dualizing, one obtains the structure of $D \kring \finj (\mathbf{a}, -)$ (cf. \cite[Corollary 6.7]{P_finmod}) with composition factors satisfying the analogous  vanishing and detection properties. 

By Example \ref{exam:indentify_Theta}, $\Theta_\mathbf{a}$ is  an isomorphism when evaluated on $\mathbf{a}$. Hence, by the above discussion, to prove the result, it suffices to show that it is not an isomorphism evaluated on $\mathbf{a+1}$. This is a straightforward calculation. 
\end{proof}

\begin{remark}
\label{rem:Feigin-Shoikhet}
The structure of $\kring \finj (\mathbf{n}, -)$, for $n \in \nat$, can also be deduced from the work of Feigin and Shoikhet \cite{MR2350124}, in which they consider certain simple representations of the Witt Lie algebra. To do this, one uses the relationship between polynomial representations of the Witt Lie algebra and $\kring \fin\op$-modules that was established  by Sam, Snowden and Tosteson in \cite{MR4788713}, where the relationship with the results of \cite{2014arXiv1406.0786W} is indicated in \cite[Section~1.7]{MR4788713}.

Now,  in \cite[Section 3.1]{MR2350124}, Feigin and Shoikhet  consider representations of the Lie algebra $W_n$ of polynomial
vector fields on an $n$-dimensional vector and its sub Lie algebra $W_n^0$ of vector fields vanishing at the origin, which is equipped with a surjection  $W_n^0 \twoheadrightarrow \mathfrak{gl}_n$.
 The proof of their classification theorem \cite[Section 3.1]{MR2350124} gives the structure of the coinduced module $\hom_{U (W_n)} (U (W_n^0), \mathfrak{gl}_n)$.  We claim that this corresponds to the structure of $\kring \finj (\mathbf{n}, -) $ via the results of \cite{MR4788713}.
\end{remark}

\begin{example}
\label{exam:case__a=0}
When $a=0$, one identifies $\kring \fs (-, \mathbf{0}) \cong \kring_\mathbf{0}$, the $\kring \fin\op$-module with value $\kring$ supported on $0$. Likewise,  $D \kring \finj (\mathbf{0}, -)$ identifies as the constant $\kring \fin\op$-module with value $\kring$. The morphism $\Theta_\mathbf{0}$ fits into the short exact sequence of $\kring \fin\op$-modules:
\[
0 
\rightarrow 
\kring_\mathbf{0}
\stackrel{\Theta_\mathbf{0}}{\rightarrow} 
\kring
\rightarrow 
\overline{\kring}
\rightarrow 
0,
\]
where $\overline{\kring}$ is $\kring$ supported on non-empty finite sets. (We stress that this is in $\kring \fin\op$-modules, as opposed to $\kring \fin$-modules; to be consistent with previous notation, these modules should be denoted $D\kring_\mathbf{0}$, $D \kring$, and $D \overline{\kring}$ respectively.) 
Now, $\overline{\kring}$ identifies as $D \Lambda^0 (\pbar)$ (by our convention) and is simple; evaluated on $\mathbf{b}$ with $b>0$, the corresponding $\kring \sym_b\op$-representation is the trivial representation $\triv_b = S_{(b)}$.
\end{example}

\section{Filtering $\kring \finj\op$-modules}
\label{sect:filter}

The purpose of this section is to introduce the natural filtration of the morphisms of $\kring \fs$ that interests us, by exploiting the $\kring \finj\op$-module structure introduced in Section \ref{sect:background}. This allows us to introduce the  wide $\kring$-linear subcategory $\primfs$ of $\kring \fs$ in Corollary \ref{cor:properties_filtration}, the main result of the section.

\subsection{Defining the filtration}

In this section, we show that a $\kring \finj\op$-module has a natural filtration. 

We start by recalling the Day convolution product on $\kring \fb\op$-modules:

\begin{definition}
\label{defn:day_convolution}
Let $\odot_{\fb\op}$ be the Day convolution product on $\kring \fb\op$-modules that is given by  $$
M \odot_{\fb\op} N : S \mapsto \bigoplus _{S_1 \amalg S_2 = S} M (S_1) \otimes N(S_2),$$
 where the sum is indexed by ordered decompositions of the finite set $S$ into two subsets $(S_1, S_2)$. 
\end{definition}

This defines a symmetric monoidal structure on $\kring \fb\op$-modules, with unit $\kring_\mathbf{0}$, the module supported on $\mathbf{0}$ with value $\kring$.

\begin{remark}
\label{rem:triv-comodules}
It is well known that the category of $\kring \finj$-modules can be described as a category of modules in the category of $\kring \fb$-modules equipped with the corresponding Day convolution product.  For example, via the Schur correspondence (as explained in \cite{MR3430359}), this description of $\kring \finj$-modules is given by \cite[Proposition 1.3.5]{MR3430359}. 

There is an analogous description of $\kring \finj\op$-modules  as a category of comodules, working in the symmetric monoidal category $(\f (\fb\op), \odot_{\fb\op} , \kring _\mathbf{0})$.  Write $\triv$ for the constant $\kring \fb\op$-module with value $\kring$. This has the structure of a bicommutative Hopf algebra in $\kring \fb\op$-modules; in particular,  $\triv$ has an underlying coaugmented coalgebra structure. One can thus consider the category of $\triv$-comodules.  Objects are $\kring \fb\op$-modules $M$ equipped with a structure morphism $M \rightarrow M \odot_{\fb\op} \triv$ that satisfies the appropriate axioms. 

The filtration of a $\kring\finj\op$-module that we construct below can be identified with the primitive filtration associated to this comodule structure.
\end{remark}

For a $\kring \finj\op$-module $M$ and a pair of finite sets $(X, Y)$, one has  the structure map $M(Y) \otimes_{\aut (Y)} \kring \finj (X, Y) \rightarrow M (X)$ and this is $\kring \aut (X)\op$-equivariant. By adjunction, this  gives the $\kring \aut (Y)\op$-equivariant morphism 
\begin{eqnarray}
\label{eqn:M_X_Y_structure_map}
M(Y) \rightarrow \hom _{\kring \aut (X)\op} (\kring \finj (X, Y), M (X)).
\end{eqnarray}

We first define the filtration at the level of sections:

\begin{definition}
\label{defn:naive_filt}
For a $\kring \finj\op$-module $M$, an integer $t \in \nat \cup \{-1\}$, and a finite set $Y$, let $M^{t} (Y) \subseteq M (Y) $ be the largest sub $\kring$-vector space such that the composite 
\[
M^{t} (Y) \subset M(Y) \stackrel{M(i)}{\rightarrow} M (X)
\]
is zero for each $(X,i)$, where $X$ is a finite set with $|Y|-|X| = t+1 $ and $i \in \finj (X, Y)$. 
\end{definition}

Since $\fin$ has a skeleton with objects indexed by their cardinality, one has the following identification:

\begin{lemma}
\label{lem:identify_M^t}
For a $\kring \finj\op$-module $M$, $t \in \nat \cup \{-1\}$, and  a finite set $Y$: 
\begin{enumerate}
\item 
\label{item:t_large}
if $|Y|\leq t$, then $M^t (Y)= M(Y)$; 
\item 
if $|Y|>t$, then 
\[
M^t (Y) = \ker \Big(
M(Y) \rightarrow \hom _{\kring \aut (X)\op} (\kring \finj (X, Y), M (X))
\Big), 
\]
where $X$ is some object with $|Y|-|X| = t+1$. 
\end{enumerate}
In particular,  $M^t (Y)$ is a sub $\kring \aut (Y)\op$-module of $M(Y)$. 
\end{lemma}

This provides the desired filtration, by the following:

\begin{proposition}
\label{prop:filter_finjop}
Let  $M$ be a $\kring \finj\op$-module.
\begin{enumerate}
\item 
For $t \in \nat \cup \{ -1 \}$, $Y \mapsto M^t (Y)$ yields a  sub $\kring\finj\op$-module of $M$. 
\item 
The subobjects $M^t$ form a natural, increasing filtration of $M$ by sub $\kring\finj\op$-modules:
\[
0 = M^{-1} \subseteq M^0 \subseteq M^1 \subseteq M^2 \subseteq \ldots \subseteq M
\]
and this filtration is exhaustive.
\end{enumerate}
\end{proposition}

\begin{proof}
It is clear from the definition that the construction of $M^t (Y) \subseteq M (Y)$ is natural with respect to $M$. 

We next establish that, if $t \leq \ell$, then $M^t (Y) \subseteq M^\ell (Y)$. Let $ i : X \hookrightarrow Y$ be any injection with $|Y|- |X| = \ell +1 $. Since $t \leq \ell$ by hypothesis, there exists a factorization of $i$ of the form $X \stackrel{i'}{\hookrightarrow} X' \stackrel{i''} {\hookrightarrow} Y$ with $|Y|- |X'|= t+1$. Thus $M(i) = M (i'') \circ M(i')$, from which one deduces that $M^t(Y)$ lies in $M^{\ell} (Y)$.

To prove that the construction $M^t$ yields a submodule of $M$, we require to show that, for any $j : U \hookrightarrow Y$ the image of $M^t (Y) \subseteq M(Y) \stackrel{j}{\rightarrow} M (U)$ lies in $M^t (U)$. If $|U|\leq t$, by Lemma \ref{lem:identify_M^t}, $M^t(U) = M(U)$ and there is nothing to show. In the case $|U|> t$, consider any $i : V \hookrightarrow U $ with $|U| - [V|=t+1$; the composite $V \stackrel{i}{\hookrightarrow} U \stackrel{j} {\hookrightarrow} Y$ has $|Y|- |V| = \ell +1$ for some $\ell \geq t$.  It follows that the composite 
 \[
 M^t (Y) 
 \subseteq M^\ell (Y) 
 \subseteq M (Y) 
 \stackrel{M (j)} {\rightarrow}
 M (U) 
 \stackrel{M(i)} {\rightarrow}
 M (V)
 \]
is zero, by the definition of $M^{\ell} (Y)$, where the first inclusion is given by the previous argument.
 This holds for any such $i$; thus, by the definition of $M^t(U)$, the image of $M^t (Y)$ in $M(U)$ lies in $M^t (U)$.

The identification of $M^{-1}$ is immediate, since any $i : X \hookrightarrow Y$ with $|Y|= |X|$ is necessarily a bijection. The fact that the filtration is exhaustive follows from the identification $M^t (Y) = M(Y)$ if $t \geq |Y| $ given in Lemma \ref{lem:identify_M^t}.
\end{proof}

There is an augmentation $\kring \finj \twoheadrightarrow \kring \fb$ that sends all non-bijective maps to zero. Restriction along the (opposite of the) augmentation gives an exact functor from $\kring \fb\op$-modules to $\kring \finj\op$-modules with essential image the full subcategory of $\kring \finj\op$-modules on which all non-bijective injections act by zero.  

We have the following interpretation of $M \mapsto M^0$ on $\kring \finj\op$-modules.

\begin{proposition}
\label{prop:primitives}
The association $M \mapsto M^0$ defines a functor from $\kring \finj\op$-modules to $\kring \fb\op$-modules. It is a right adjoint to the above restriction functor.

The adjunction counit identifies as the natural inclusion  $M^0 \hookrightarrow M$ of $\kring \finj\op$-modules. 
\end{proposition} 

\begin{proof}
By definition, for a $\kring \finj\op$-module $M$, $M^0$ is the largest submodule of $M$ with the property that $M^0(i)$ is zero for every non-bijective injection $i$. Thus $M \mapsto M^0$ induces a functor from $\kring \finj\op$-modules to $\kring \fb\op$. Moreover, it is immediate that this is right adjoint to the restriction functor and that the adjunction counit corresponds to the canonical inclusion $M^0 \hookrightarrow M$. 
\end{proof}

\begin{remark}
\label{rem:primitives}
The submodule $M^0$ corresponds to the primitives of 
$M$  when considered as a comodule in $\kring \fb\op$-modules (see Remark \ref{rem:triv-comodules}).
\end{remark}

\subsection{Relating terms of the filtration}
In this section, we use the $\kring \finj\op$-module structure of $M$ to provide further relations between the terms $M^t$ of its filtration. 
For this, we introduce the following:

\begin{notation}
\label{nota:finj_deg_n}
For $n \in \nat$, set 
$$
\finj\dg{n} (X, Y):= 
\left\{ 
\begin{array}{ll}
\finj (X,Y) & |Y| - |X|= n \\
\emptyset & \mbox{otherwise.}
\end{array}
\right.
$$
\end{notation}

\begin{remark}
Forgetting structure, we may consider $\kring \finj$ as a $\kring \fb$-bimodule. By construction 
$\kring \finj\dg{n}$ is a sub $\kring \fb$-bimodule and 
there is an isomorphism of $\kring \fb$-bimodules $
\kring \finj \cong \bigoplus _{n \in \nat} \kring \finj\dg{n}$. Observe that $\kring \finj\dg{0}$ identifies as $\kring \fb$.
 This provides an $\nat$-grading of the morphisms of the $\kring$-linear category $\kring \finj$.
\end{remark}

We can form $\hom_{\kring \fb \op} ( \kring \finj\dg{n}, M)$ for any $\kring \fb\op$-module $M$. The following is immediate:

\begin{lemma}
\label{lem:induction}
For $M$ a $\kring \fb\op$-module and $n \in \nat$, the underlying $\kring \fb$-module structure of $\kring \finj\dg{n}$ induces a natural $\kring\fb\op$-module structure on $ \hom_{\kring \fb \op} ( \kring \finj\dg{n}, M )$.

Evaluated on $Y$, if $|Y|\leq n$, this is zero; otherwise it  is isomorphic to 
$$\hom _{\kring \aut (X)\op} (\kring \finj (X, Y), M (X)),$$
where $X$ is a set of cardinality $|Y|-n$.
\end{lemma}

We identify $ \hom_{\kring \fb \op} ( \kring \finj\dg{n}, - )$ using the 
Day 
convolution product $\odot_{\fb\op}$ that was recalled in Definition \ref{defn:day_convolution}.

For $n \in \nat$, we consider the trivial $\kring \sym_n$-module $\triv_n$  as a $\kring \fb\op$-module supported on sets of cardinality $n$. Hence, for any $\kring \fb\op$-module $M$ one can form the $\kring \fb\op$-module $M \odot_{\fb\op} \triv_n$. 

\begin{lemma}
\label{lem:identify_induction}
For $M$ a $\kring \fb\op$-module, there is a natural isomorphism of $\kring \fb\op$-modules
\[
\hom_{\kring \fb \op} ( \kring \finj\dg{n}, M )
\cong 
M \odot_{\fb\op} \triv_n.
\]
\end{lemma}

\begin{proof}
If $b <n$, then $\hom_{\kring \fb \op} ( \kring \finj \dg{n} (-,\mathbf{b}), M )$ and $(M \odot_{\fb\op} \triv_n)(\mathbf{b})$ are clearly zero. It remains to treat the case $b \geq n$;  set $a:= b - n$. 

By Lemma \ref{lem:induction}, there is an isomorphism of $\kring \sym_b\op$-modules:
\[
\hom_{\kring \fb \op} ( \kring \finj\dg{n} (-,\mathbf{b}), M )
\cong 
\hom_{\kring \sym_a\op} ( \kring \finj  (\mathbf{a},\mathbf{b}), M(\mathbf{a}) ).
\]

Consider the Young subgroup $\sym_a \times \sym_{n} \subset \sym_b$ corresponding to the canonical inclusion $\mathbf{a}\subset \mathbf{b}$, using that $\mathbf{b} \backslash \mathbf{a} \cong \mathbf{n}$. The $\kring \sym_a\op \otimes \kring \sym_b$-module $\kring \finj  (\mathbf{a},\mathbf{b})$  identifies as 
the permutation module
$
\kring (\sym_b /\sym_n)$, 
 where $\kring \sym_a$ acts on the right via the inclusion $\sym_a \subset \sym_b$. From this one deduces the isomorphism:
\[
\hom_{\kring \sym_a\op} ( \kring \finj  (\mathbf{a},\mathbf{b}), M(\mathbf{a}) )
\cong 
\hom_{\kring (\sym_a \times \sym_n)\op} (\kring \sym_b, M (\mathbf{a}) \boxtimes \triv_n)).
\]
Since, for finite groups, induction is isomorphic to coinduction, this provides the natural isomorphism:
\[
\hom_{\kring \sym_a\op} ( \kring \finj  (\mathbf{a},\mathbf{b}), M(\mathbf{a}) )
\cong 
\big(M (\mathbf{a}) \boxtimes \triv_n\big) \otimes_{\kring (\sym_a \times \sym_n)} \kring \sym_b.
\]
The right hand side is naturally isomorphic to $(M \odot_{\fb\op} \triv_n) (\mathbf{b})$, as required.
\end{proof}

\begin{remark}
\label{rem:towards_comodule_structure}
Suppose that $M$ is a $\kring \finj\op$-module. Fixing the difference of the cardinalities $|Y|-|X|=n$, for some $ n \in \nat$, the structure maps (\ref{eqn:M_X_Y_structure_map}) assemble to define a natural transformation of $\kring \fb\op$-modules
\begin{eqnarray}
\label{eqn:structure_map_M_n}
M \rightarrow \hom_{\kring \fb \op} ( \kring \finj\dg{n}, M )
\end{eqnarray}
using the restricted $\kring \fb\op$-module structure of $M$ in the codomain. 

By Lemma \ref{lem:identify_induction}, this morphism can be rewritten as the morphism of $\kring \fb\op$-modules 
$
M \rightarrow M \odot_{\fb\op} \triv _n 
$. 
Since $\triv \cong \bigoplus_{n \in \nat} \triv_n$ as $\kring \fb\op$-modules, these assemble to the structure morphism
$$
M
\rightarrow 
M \odot_{\fb\op} \triv. 
$$
This is the structure morphism of the $\triv$-comodule structure of Remark \ref{rem:triv-comodules}. The coassociativity and counit properties of this follow from the associativity and unit properties of the $\kring \finj\op$-module structure.  This leads to the equivalence  of Remark \ref{rem:triv-comodules}.
\end{remark}

If $M$ is a $\kring \finj\op$-module, then, by Lemma \ref{lem:identify_M^t}, for $t \in \nat \cup \{-1 \}$, the underlying $\kring \fb\op$-module of $M^t$ is the kernel of the structure map (\ref{eqn:structure_map_M_n}) with $n=t+1$:
\[
M \rightarrow \hom_{\kring \fb \op} ( \kring \finj \dg{t+1}, M ).
\]

Then Proposition \ref{prop:filter_finjop} refines to:

\begin{proposition}
\label{prop:refined_filtration_result}
For $M$ a $\kring \finj\op$-module and $n \in \nat$, the structure map (\ref{eqn:structure_map_M_n})  
restricts for $t \in \nat$ to the natural map 
$$
M^t
\rightarrow 
\hom_{\kring \fb \op} ( \kring \finj \dg{n}, M^{t-n} ), 
$$
where $M^{t-n}$ is taken to be zero if $t-n <0$.

In particular, this yields the natural injection of $\kring \fb\op$-modules:
\begin{eqnarray}
\label{eqn:Mn/n-1_to_M0}
M^n/M^{n-1} 
\hookrightarrow 
\hom_{\kring \fb \op} ( \kring \finj \dg{n}, M^{0} )
\cong 
M^0 \odot_{\fb\op} \triv_n.
\end{eqnarray}
\end{proposition}

\begin{proof}
The first statement is proved by refining the proof of Proposition \ref{prop:filter_finjop}, taking into account the cardinalities.

Since $M^{n-1}$ is the kernel of $M
\rightarrow 
\hom_{\kring \fb \op} ( \kring \finj \dg{n}, M )$, the transformation (\ref{eqn:structure_map_M_n}) factors across the injection
\[
M/M^{n-1} \hookrightarrow \hom_{\kring \fb \op} ( \kring \finj \dg{n}, M).
\]
Upon restriction to $M^n/M^{n-1}$ this yields the required injection; this takes values in $\hom_{\kring \fb \op} ( \kring \finj \dg{n}, M^0)$, by the first statement.

Naturality with respect to $M$ is clear. 
\end{proof}

\subsection{Filtering  $\kring \fs$}

By Proposition \ref{prop:quotient_fs_finj}, for a fixed finite set $X$, $\kring \fs (-,X)$ has a $\kring \aut (X) \otimes \kring \fin\op$-module structure. Restricting along $\kring \finj \subset \kring \fin$, it follows that $\kring \fs (-, X)$ is a $\kring \aut (X) \otimes \kring \finj\op$-module. The structure maps of this module, in the sense of the previous subsection, are given by Proposition \ref{prop:adjoint_XY}; these are used in defining the filtration in Corollary \ref{cor:filt_qfs} below.

Proposition \ref{prop:filter_finjop} together with Proposition \ref{prop:refined_filtration_result} 
yield:

\begin{corollary}
\label{cor:filt_qfs}
For a finite set $X$, 
the filtration of Proposition \ref{prop:filter_finjop} gives an increasing, exhaustive filtration by sub $\kring \aut (X) \otimes \kring \finj\op$-modules:
\[
0=\kring \fs (-, X)^{-1} \subseteq \kring \fs (-, X)^0  \subseteq \kring \fs (-, X)^1  \subseteq\kring \fs (-, X )^2  \subseteq\ldots \subseteq \kring \fs (-, X).
\]

Moreover, the $\kring \finj\op$-module structure induces an injection of 
\begin{eqnarray*}
\kring \fs (-,X)^t / \kring\fs (-, X)^{t-1}
&\hookrightarrow&
\hom_{\kring \fb \op} (\kring \finj\dg{t} (-,-) , \kring \fs (-,X)^0)
\\
&\phantom{\hookrightarrow}&
\cong
\kring\fs (-, X)^0 \odot_{\fb\op} \triv_t.
\end{eqnarray*}
\end{corollary}

We introduce the following notation (which is extended and refined in Notation \ref{nota:primfs_ell}): 

\begin{notation}
\label{nota:primfs}
Write $\primfs$ for $\kring \fs (-,-)^0$, so that $\primfs (X,Y) = \kring \fs (-,Y)^0 (X)$ for finite sets $X,Y$. 
\end{notation}

\begin{remark}
\label{rem:primfs_explicit}
The definition of $\primfs$ can be made explicit as follows. For $a, b \in \nat$, we have 
\begin{eqnarray}
\label{eqn:primfs_explicit}
\primfs (\mathbf{b}, \mathbf{a}) :=
\ker \big( 
\kring \fs (\mathbf{b} , \mathbf{a}) 
\rightarrow 
\bigoplus_{\substack{X \subset \mathbf{b} \\ |X|= b-1 }} \kring \fs (X, \mathbf{a}) 
\big),
\end{eqnarray}
where the map is defined by the $\kring \finj \op$-module structure of $\kring \fs (-, \mathbf{a})$. 

This coincides with $H^0 (\mathbb{C} (\mathbf{b}, \mathbf{a}))$ that is studied in  \cite[Section 6]{MR4518761}, where $\mathbb{C}(\mathbf{b}, \mathbf{a})$ is defined in terms of  the Koszul complex associated to the $\kring \finj\op$-module $\kring \fs (-, \mathbf{a})$. 
\end{remark}

Forgetting structure, $\kring \fs$ can be considered as a $\kring \fb$-bimodule and $\primfs$ is a sub $\kring \fb$-bimodule. The following is immediate:

\begin{lemma}
\label{lem:prim_contains_kfb}
The inclusion $\kring \fb \hookrightarrow \kring \fs$ factorizes across $\primfs$ to give the inclusion  $\kring \fb \hookrightarrow \primfs$ of $\kring \fb$-bimodules. This is an isomorphism when evaluated on $(\mathbf{n},\mathbf{n})$, for any $n \in \nat$.
\end{lemma}

\begin{remark}
\label{rem:primfs_non_trivial}
The structure of the $\kring \fb\op$-bimodule $\primfs$ is much richer than that of $\kring \fb$: the above inclusion is far from being a bijection, as is shown by the results of \cite[Section 6]{MR4518761} (and can also be deduced from Theorem \ref{thm:primfs_bimodule} below).

As stated in Remark \ref{rem:primfs_explicit}, $\primfs (\mathbf{b}, \mathbf{a})$ identifies with $H^0 (\mathbb{C} (\mathbf{b}, \mathbf{a}))$ of  \cite[Section 6]{MR4518761}. The latter is shown to be highly non-trivial by \cite[Proposition 6.18]{MR4518761}, which gives a recursive expression for the dimension of this vector space. Moreover, for positive integers $b>a$,  \cite[Theorem 6.20]{MR4518761} identifies the isomorphism class of this $\kring (\sym_b\op \times \sym_a)$-module.
\end{remark}

\begin{example}
\label{exam:primfs}
The following results are contained in \cite[Section 6]{MR4518761}.
\begin{enumerate}
\item 
For $a=1$, $\primfs (\mathbf{b}, \mathbf{1})$ is isomorphic to $\kring$  for $b=1$ and is zero for $b>1$. 
\item 
For $a=2$, $\primfs (\mathbf{b}, \mathbf{2})$ is isomorphic (as a $\kring (\sym_b\op \times \sym_2)$-module) to $\kring \sym_2$ (with the regular bimodule structure) for $b=2$ and, for $b>2$,  to 
$\triv_b \boxtimes \sgn_2$ if $b$ is odd and $\triv_b \boxtimes \triv_2$ if $b$ is even. 

This can be seen as follows.  If $b\geq 2$, a surjection $f : \mathbf{b} \twoheadrightarrow \mathbf{2}$ is uniquely determined by the fibre $f^{-1} (1)$, which must be a subset $\emptyset \subsetneq U \subsetneq \mathbf{b}$; write $f_U$ for the corresponding surjection. Then, for $b>2$,   $\primfs (\mathbf{b}, \mathbf{2})$ is one-dimensional, and a generator is given by 
\[
\sum_{\emptyset \subsetneq U \subsetneq \mathbf{b}} (-1)^{|U|} [f_U].
\] 
Indeed, it is straightforward to show that the vector space is at most one dimensional and that the exhibited element lies in the kernel of (\ref{eqn:primfs_explicit}). The identification of the representation also follows from this.
\end{enumerate}
\end{example}

The  following result is the key input to Corollary \ref{cor:properties_filtration}:

\begin{proposition}
\label{prop:fundamental_primfs}
For finite sets $X$, $Y$, composition in $\kring \fs$ induces a morphism of $\kring \finj\op$-modules:
\[
\primfs (X,Y) \otimes \kring \fs (-, X) \rightarrow \kring \fs(- , Y). 
\]
For $\ell \in \nat$, this restricts to a morphism of $\kring \finj\op$-modules:
\[
\primfs (X,Y) \otimes \kring \fs (-, X)^\ell \rightarrow \kring \fs(-, Y)^\ell.
\]
\end{proposition}

\begin{proof}
Recall that the $\kring \finj\op$-module structure on $\kring \fs (-, X)$ (respectively with $Y$ in place for $X$) is the restriction of the  $\kring \fin\op$-module structure provided by Proposition \ref{prop:quotient_fs_finj}, so that the quotient map $\kring \fin (-, X) \twoheadrightarrow \kring \fs (-, X)$ is a morphism of $\kring \fin\op$-modules and hence, by restriction, of $\kring \finj\op$-modules (resp. for $Y$). 

Now, using the inclusions $\kring \fs \hookrightarrow \kring \fin$ and $\kring \finj \hookrightarrow \kring \fin$ of wide $\kring$-linear subcategories, for any finite sets $A$, $B$, associativity of composition in $\kring \fin$ yields the following commutative diagram:
\[
\xymatrix{
\kring \fs (X, Y) \otimes \kring \fs (B, X) \otimes \kring \finj (A, B) 
\ar[r]
\ar[d]
&
\kring \fs (B, Y) \otimes \kring \finj (A, B) 
\ar[d]
\\
\kring \fs (X,Y) \otimes \kring \fin (A, X) 
\ar[r]
&
\kring \fin (A,Y), 
}
\]
In which the maps are given by the appropriate compositions, passing to the larger category $\kring \fin$ where required. 

The subtlety comes from the fact that the `action' of $\kring \fs(X,Y)$ is not in general compatible with the projections $\kring \fin (-,X)\twoheadrightarrow \kring \fs (-, X)$ and $\kring \fin (-,Y)\twoheadrightarrow \kring \fs (-, Y)$. Namely, in general, the following diagram does  not commute:
\[
\xymatrix{
\kring \fs (X,Y) \otimes \kring \fin (A, X) 
\ar[r]
\ar@{}[rd]|{ \not\circlearrowleft}
\ar@{->>}[d]
&
\kring \fin (A,Y)
\ar@{->>}[d]
\\
\kring \fs (X,Y) \otimes \kring \fs (A, X)
\ar[r]
&
\kring \fs (A,Y), 
}
\]
since a non-surjective map $A \rightarrow X$ may become surjective after composing with a surjection $X \twoheadrightarrow Y$. 

Restricting to $\primfs (X, Y)\subset \kring \fs (X,Y)$ resolves this issue. Namely, we claim that the restricted diagram 
\[
\xymatrix{
\primfs (X,Y) \otimes \kring \fin (A, X) 
\ar[r]
\ar@{->>}[d]
&
\kring \fin (A,Y)
\ar@{->>}[d]
\\
\primfs (X,Y) \otimes \kring \fs (A, X)
\ar[r]
&
\kring \fs (A,Y)
}
\]
always commutes. 

This  is seen as follows. Write $\overline{\kring \finj}$ for the kernel of the augmentation $\kring \finj \twoheadrightarrow \kring \fb$ that is induced by sending the non-bijective maps to zero. One checks  that  
the map induced by composition in $\kring \fin$
\[
\overline{\kring \finj} (-, X) \otimes_\fb \kring \fin (A, -) \rightarrow \kring \fin (A, X) 
\]
surjects to the $\kring$-linear span of the non-surjective maps from $A$ to $X$. 
 
Now, for any finite set $U$,  the map induced by composition in $\kring \fin$ 
\[
\primfs ( X, Y) \otimes \overline{\kring \finj} (U, X) \rightarrow \kring  \fin(U, Y) 
\]
maps to the $\kring$-linear span of the non-surjective maps from $U$ to $Y$, since the composite with $\kring \fin (U, Y) \twoheadrightarrow \kring \fs (U, Y)$ is zero, by definition of $\primfs$. The claimed commutativity follows by putting these ingredients together, using associativity of composition in $\kring \fin$.

Thus, after restricting to $\primfs$, one can paste the two commutative diagrams together to obtain the commutative diagram
\begin{eqnarray}
\label{eqn:kfs_kfinjop-module}
\xymatrix{
\primfs (X, Y) \otimes \kring \fs (B, X) \otimes \kring \finj (A, B) 
\ar[r]
\ar[d]
&
\kring \fs (B, Y) \otimes \kring \finj (A, B) 
\ar[d]
\\
\primfs (X,Y) \otimes \kring \fs (A, X) 
\ar[r]
&
\kring \fs (A,Y), 
}
\end{eqnarray}
in which the vertical maps are given by the $\kring \finj\op$-structure maps and the horizontal maps correspond to the action of $\primfs (X,Y)$. This establishes the first statement. 

The second statement then follows  from the naturality of the filtration of $\kring \finj\op$-modules, as stated in Proposition \ref{prop:filter_finjop}.  
\end{proof}

This gives the following fundamental properties of the filtration of morphisms of $\kring \fs$, in particular, highlighting the structure of $\primfs$.

\begin{corollary}
\label{cor:properties_filtration}
\ 
\begin{enumerate}
\item 
The $\kring$-linear category $\primfs$ is  a wide $\kring$-linear subcategory of $\kring \fs$ that contains $\kring \fb$ as a wide $\kring$-linear subcategory;
\item 
$\kring \fs $ is  a $\primfs \otimes \kring \finj\op$-module.
\item 
The filtration $\kring \fs (-, -) ^\ell$ is a filtration by $\primfs \otimes \kring \finj\op$-modules. 
\item 
The induced inclusion 
\[
\kring \fs (-,-)^{\ell} / \kring \fs (-,-)^{\ell -1} 
\hookrightarrow 
\hom_{\kring \fb\op} (\kring \finj\dg{\ell} (-,-), \primfs) 
\]
is a morphism of $\primfs \otimes \kring \fb\op$-modules. 
\end{enumerate}
\end{corollary}

\begin{proof}
The first point follows by taking $\ell=0$ in the second statement of Proposition \ref{prop:fundamental_primfs}. 
 Then, knowing that $\primfs$ is a $\kring$-linear subcategory of $\kring \fs$, the commutative diagram (\ref{eqn:kfs_kfinjop-module}) implies that $\kring \fs $ is  a $\primfs \otimes \kring \finj\op$-module, using associativity of composition.

The remaining points follow from the naturality of the filtration of Proposition \ref{prop:filter_finjop}. 
\end{proof}

\section{The subquotients of the filtration of $\kring \fs$}
\label{sect:subquotients}

The purpose of this section is to identify the subquotients $\kring \fs (-,-)^\ell/ \kring \fs (-,-)^{\ell -1}$ of the filtration introduced in Section \ref{sect:filter} as $\primfs \otimes \kring \fb\op$-modules. 

The main result is Theorem \ref{thm:coker_theta_ell}, which determines these subquotients in terms of $\primfs$. This then leads to Corollary \ref{cor:kring_fs}, which gives the explicit relationship between $\primfs$ and $\kring \fs$ that is exploited in Section \ref{sect:primfs} to calculate the underlying $\kring \fb$-bimodule of $\primfs$.

\subsection{Preliminaries}
To simplify the notation, we introduce the following. 

\begin{notation}
\label{nota:primfs_ell}
For $\ell \in \nat$, denote by
\begin{enumerate}
\item 
$\kring \fs^\ell$  the $\primfs \otimes \kring \fb\op$-module $\kring \fs (-,-)^\ell$ (by convention, set  $\kring \fs^{-1} :=0$); 
\item 
$\kring \fs^{\ell/(\ell-1)}$  the quotient $\kring \fs ^\ell/ \kring \fs ^{\ell -1}$; 
\item 
$\theta^\ell$ the morphism $\kring \fs^{\ell/(\ell-1)} 
\hookrightarrow 
\hom_{\kring \fb\op} (\kring \finj\dg{\ell}, \primfs) $ of Corollary \ref{cor:properties_filtration}.
\end{enumerate}
\end{notation}

For $\ell, b\in \nat$, from the definition of $\kring \finj\dg{\ell}$, one has that $\kring \finj\dg{\ell}(- , \mathbf{b})$ is zero if $\ell >b$ and otherwise is supported on $\mathbf{c}$, where $c = b- \ell$.  Thus, evaluating $\theta^\ell$ on $(\mathbf{b}, \mathbf{a})$, for $a,b \in \nat$, gives 
\[
\kring \fs^{\ell/(\ell-1)} (\mathbf{b}, \mathbf{a}) 
\hookrightarrow 
\hom_{\kring \sym_c\op} (\kring \finj\dg{\ell} (\mathbf{c}, \mathbf{b}) , \primfs (\mathbf{c}, \mathbf{a}))
\]
where $c = b - \ell$. 

We first check the behaviour of the domain and codomain of $\theta^\ell (\mathbf{b}, \mathbf{a})$ for $\ell \geq b-a$:

\begin{lemma}
\label{lem:vanishing_identification}
For $\ell, a \leq  b \in \nat$, 
\begin{enumerate}
\item 
if $\ell = b-a$, then, as $\kring \sym_a \otimes \kring \sym_b\op$-modules,
\[
    \kring \fs^{b-a} (\mathbf{b}, \mathbf{a}) = \kring \fs (\mathbf{b}, \mathbf{a})
\]
and
\[
    \hom_{\kring \fb\op} (\kring \finj\dg{b-a}, \primfs)(\mathbf{b}, \mathbf{a}) \cong \hom_\kring (\kring \finj (\mathbf{a}, \mathbf{b}), \kring).
\]
\item 
if $\ell > b-a$, then  $\kring \fs^{\ell/(\ell-1)} (\mathbf{b}, \mathbf{a}) $ and $\hom_{\kring \fb\op} (\kring \finj\dg{\ell}, \primfs)(\mathbf{b}, \mathbf{a})$ are both zero.
\end{enumerate} 
\end{lemma}

\begin{proof}
For any $n \in \nat$, $\kring \fs (\mathbf{n}, \mathbf{a})$ vanishes if $n < a$. This  implies the vanishing of $\kring \fs^{\ell} (\mathbf{n}, \mathbf{a})$, $\kring \fs^{\ell/(\ell-1)} (\mathbf{n}, \mathbf{a})$ and of $\primfs (\mathbf{n}, \mathbf{a}$) for any $n$.

The equality  $\kring \fs^{b-a} (\mathbf{b}, \mathbf{a}) = \kring \fs (\mathbf{b}, \mathbf{a})$ follows directly  from these observations (compare Lemma \ref{lem:identify_M^t} (\ref{item:t_large})). For the second isomorphism, using the identification of the support of $\kring \finj\dg{b-a}(-, \mathbf{b})$ given before the statement of the Lemma, we have
$$
\hom_{\kring \fb\op} (\kring \finj\dg{b-a}, \primfs)(\mathbf{b}, \mathbf{a}) 
\cong 
\hom_{\kring \sym_a\op} (\kring \finj (\mathbf{a},\mathbf{b}), \primfs (\mathbf{a}, \mathbf{a})) 
$$ 
and $\primfs (\mathbf{a}, \mathbf{a}) = \kring \fs (\mathbf{a}, \mathbf{a}) \cong \kring \sym_a$, by Lemma \ref{lem:prim_contains_kfb}. The isomorphism in the statement then follows, as in Corollary \ref{cor:adjoint_X} (see the argument given in the proof of that result).

The equality  $\kring \fs^{b-a} (\mathbf{b}, \mathbf{a}) = \kring \fs (\mathbf{b}, \mathbf{a})$  implies  that $\kring \fs^{\ell/(\ell-1)} (\mathbf{b}, \mathbf{a})=0$ if $\ell > b-a$. Finally,  $\hom_{\kring \fb\op} (\kring \finj\dg{\ell}, \primfs)(\mathbf{b}, \mathbf{a})$ identifies as 
$\hom_{\kring \sym_c\op} (\kring \finj\dg{\ell} (\mathbf{c}, \mathbf{b}) , \primfs (\mathbf{c}, \mathbf{a}))$ with $c= b - \ell$. If $\ell > b-a$, then $c <a$, hence $\primfs (\mathbf{c}, \mathbf{a})=0$, which gives 
$$
 \hom_{\kring \fb\op} (\kring \finj\dg{\ell}, \primfs)(\mathbf{b}, \mathbf{a})=0,
$$ as required.
\end{proof}

For $a \leq b \in \nat$, taking $\ell = b-a$, we may therefore form the composite 
\begin{eqnarray}
\label{eqn:composite_theta_l=b-a}
\kring \fs (\mathbf{b}, \mathbf{a}) 
= 
\kring \fs^{b-a} (\mathbf{b}, \mathbf{a}) 
\twoheadrightarrow 
\kring \fs^{(b-a) / (b-a-1)}(\mathbf{b}, \mathbf{a}) 
\stackrel{\theta^{b-a} } {\longrightarrow}
 \hom_\kring (\kring \finj (\mathbf{a}, \mathbf{b}), \kring),
\end{eqnarray}
using the identification from Lemma \ref{lem:vanishing_identification} for the codomain.

\begin{proposition}
\label{prop:identify_composite}
The composite (\ref{eqn:composite_theta_l=b-a}) identifies with the natural transformation 
\[
\Theta_\mathbf{a} : 
\kring \fs (-, \mathbf{a}) \rightarrow D \kring \finj (\mathbf{a}, - ) 
\]
of Corollary \ref{cor:adjoint_X} evaluated on $\mathbf{b}$. 
\end{proposition}

\begin{proof}
This is essentially tautological: both maps are defined by the same construction. 
\end{proof}

\subsection{The main result}
Recall that $\theta^\ell$ is the morphism of $\primfs \otimes \kring \fb\op$-modules 
$$\kring \fs^{\ell/(\ell-1)} 
\hookrightarrow 
\hom_{\kring \fb\op} (\kring \finj\dg{\ell}, \primfs) .$$

Using Lemma \ref{lem:identify_induction}, the codomain identifies as 
$
\primfs \odot_{\fb\op} \triv_\ell
$. This has the structure of a $\primfs \otimes \kring \fb\op$-module: evaluating on $\mathbf{a}$ (considered as an object of $\primfs$) gives the $\kring \fb\op$-module 
$
\primfs (-, \mathbf{a}) \odot_{\fb\op} \triv_\ell
$; 
a morphism of $\primfs (\mathbf{a}, \mathbf{a}')$ induces a  morphism of $\kring \fb\op$-modules
$$
\primfs (-, \mathbf{a}) \odot_{\fb\op} \triv_\ell
\rightarrow 
\primfs (-, \mathbf{a}') \odot_{\fb\op} \triv_\ell
$$
 by naturality of the convolution product.

Using this structure, $\theta^\ell$ identifies as a monomorphism of 
$\primfs \otimes \kring \fb\op$-modules
\[
\theta^{\ell}: 
\kring \fs^{\ell/(\ell-1)} 
\hookrightarrow 
\primfs \odot_{\fb\op} \triv_\ell.
\]

\begin{theorem}
\label{thm:coker_theta_ell}
For $\ell >0$, there is a short exact sequence of $\primfs \otimes \kring \fb\op$-modules
\[
0
\rightarrow 
\kring \fs^{\ell/(\ell -1)}
\stackrel{\theta^\ell}{\longrightarrow} 
\primfs \odot_{\fb\op} \triv_\ell
\rightarrow 
\bigoplus_{a\in \nat}  \sgn_a \boxtimes S_{(\ell, 1^a)} 
\rightarrow 
0,
\]
where the cokernel is a $\primfs \otimes \kring \fb\op$-module by restriction of its $\kring \fb \otimes \kring \fb\op$-module structure using the augmentation $\primfs \subset \kring \fs \twoheadrightarrow \kring \fb$.
\end{theorem}

\begin{proof}
To identify the underlying $\kring \fb$-bimodule of the cokernel of $\theta^\ell$, 
it suffices to consider behaviour when evaluated on $\mathbf{b} \in \ob \kring \fb\op$, for each $b \in \nat$.

For fixed $b$, one has the morphism $\theta^\ell (\mathbf{b}, -) $ of $\primfs \otimes \kring \sym_b\op$-modules
\[
\kring \fs^{\ell / (\ell -1)} (\mathbf{b}, -) \rightarrow (\primfs \odot_{\fb\op} \triv_\ell) (\mathbf{b}, -).
\]
Evaluated on $\mathbf{d} \in \ob \primfs$, by Lemma \ref{lem:vanishing_identification}, for $d > b - \ell$, both terms are zero and this is an isomorphism. Thus, there can only be a non-zero contribution to the cokernel if $b - \ell \geq 0$, which we suppose henceforth. 

Set $a:= b - \ell$. Evaluated on $(\mathbf{b}, \mathbf{a})$, by Proposition \ref{prop:identify_composite} the cokernel of $\theta^\ell (\mathbf{b}, \mathbf{a})$ is isomorphic to the cokernel of $\Theta_\mathbf{a} (\mathbf{b})$. By Theorem \ref{thm:coker_Theta} together with Lemma \ref{lem:identify_Lambda^a_pbar}, this is isomorphic to $\sgn_a \boxtimes S_{(\ell, 1^a)} $ as a $\kring \sym_a  \otimes \kring \sym_b\op$-module. As $b$ varies, this will account for the entire cokernel of $\theta^\ell$, as shown below. 

Namely, to conclude  it suffices to show that $\theta^\ell (\mathbf{b}, \mathbf{d})$ is an isomorphism if $d < b - \ell$ (i.e., $d<a$ for $a$ as given). This is proved by using that $\theta^\ell$ is a morphism of $\primfs$-modules, as explained below.

One has the isomorphism of $\primfs$-modules
\[
(\primfs \odot_{\fb\op} \triv_\ell) (\mathbf{b}, -) \cong \primfs (\mathbf{a}, -) \otimes_{\sym_a} D \kring \finj (\mathbf{a}, \mathbf{b}).
\]
The right hand side (and hence the left) is clearly generated as a $\primfs$-module by its values on $\mathbf{a}$.

Since $\kring$ has characteristic zero, the inclusion $\mathrm{Image} (\Theta_{\mathbf{a}} (\mathbf {b})) \hookrightarrow D \kring \finj (\mathbf{a}, \mathbf{b})$ induces the inclusion of $\primfs \otimes \kring \sym_b\op$-modules:
\begin{eqnarray}
\label{eqn:inclusion_image_Theta}
\primfs (\mathbf{a}, -) \otimes_{\sym_a} \mathrm{Image} (\Theta_{\mathbf{a}} (\mathbf {b}))
\hookrightarrow 
\primfs (\mathbf{a}, -) \otimes_{\sym_a} D \kring \finj (\mathbf{a}, \mathbf{b}).
\end{eqnarray}
Since $\theta^{\ell}$ is a morphism of $\primfs$-modules, it follows that $\primfs (\mathbf{a}, -) \otimes_{\sym_a} \mathrm{Image} (\Theta_{\mathbf{a}} (\mathbf {b}))$ lies in the image of $\theta^{\ell}(\mathbf{b}, -)$, since this is true when evaluated on $\mathbf{a}$. Thus,  it suffices to show that (\ref{eqn:inclusion_image_Theta}) is an isomorphism when evaluated on $\mathbf{d}$ for $d < a$. 

By Theorem \ref{thm:coker_Theta} and Lemma \ref{lem:identify_Lambda^a_pbar} (as above), the morphism (\ref{eqn:inclusion_image_Theta}) has cokernel 
\[
\primfs (\mathbf{a}, -) \otimes_{\sym_a} ( \sgn_a \boxtimes S_{(\ell, 1^a) } )
\cong 
(\primfs (\mathbf{a}, -) \otimes_{\sym_a} \sgn_a) \boxtimes S_{(\ell, 1^a) }
\]
considered as a $\primfs \otimes \kring \sym_b\op$-module.

Now, $\primfs \hookrightarrow \kring \fs$ induces $\primfs (\mathbf{a}, -) \otimes_{\sym_a} \sgn_a \hookrightarrow \kring \fs (\mathbf{a}, -) \otimes_{\sym_a} \sgn_a$. A standard and straightforward calculation shows that, for $d < a$, one has 
$
 \kring \fs (\mathbf{a}, \mathbf{d}) \otimes_{\sym_a} \sgn_a = 0$.
 This implies that the cokernel of (\ref{eqn:inclusion_image_Theta}) is zero when evaluated on $\mathbf{d}$ for $d<a$, hence (\ref{eqn:inclusion_image_Theta}) is an isomorphism for such $d$, as required. 
 
The above shows that the underlying $\kring \fb$-bimodule of the cokernel is isomorphic to  $\bigoplus_{a\in \nat}  \sgn_a \boxtimes S_{(\ell, 1^a)} $, which arises as $\bigoplus_{a \in \nat} \mathrm{coker}\  \Theta_{\mathbf{a}} (\mathbf{b})$, where $b= a+\ell$.
 
It remains to shows that the $\primfs$-module structure arises from the $\kring \fb$-module structure by restriction along the augmentation $\primfs \rightarrow \kring \fb$; equivalently, that the augmentation ideal of $\primfs$ acts trivially.

For this, consider the contribution from $a \in \nat$ and the action of a morphism of $\primfs (\mathbf{a}, \mathbf{d})$ for $d<a$. Take $b:= a + \ell$; naturality of $\theta^\ell (\mathbf{b}, -)$ gives the commutative diagram:
\[
\xymatrix{
\kring \fs^{\ell / (\ell -1)} (\mathbf{b}, \mathbf{a})
\ar[rr]^(.45){\theta^\ell (\mathbf{b}, \mathbf{a})}
\ar[d]
&& 
 (\primfs \odot_{\fb\op} \triv_\ell) (\mathbf{b}, \mathbf{a})
 \ar[d]
 \\
\kring \fs^{\ell / (\ell -1)} (\mathbf{b}, \mathbf{d})
\ar[rr]_(.45){\theta^\ell (\mathbf{b}, \mathbf{d})}
& &
 (\primfs \odot_{\fb\op} \triv_\ell) (\mathbf{b}, \mathbf{d}),
 }
 \]
 where the vertical maps are given by the $\primfs$-naturality. As established above, the lower horizontal map is an isomorphism. Hence the induced map on the cokernel is zero. This shows that the augmentation ideal acts trivially, as required.
\end{proof}

\begin{remark}
\label{rem:theorem_coker_theta_ell}
The case $\ell =0$  corresponds to the isomorphism $\primfs \cong \primfs \odot_{\fb\op} \triv_0$. It has  not been included in the statement of Theorem \ref{thm:coker_theta_ell}  to avoid complicating the proof. To treat the case $\ell=0$ requires modifying the argument to include the case  $b=a$, for which $\Lambda^a (\pbar) (\mathbf{b})=0$ by Lemma \ref{lem:identify_Lambda^a_pbar}.
\end{remark}

If one forgets structure, by restricting from $\primfs \otimes \kring \fb\op$-modules to $\kring \fb$-bimodules, one has the isomorphism 
$
\kring \fs
\cong 
\bigoplus_{\ell \in \nat}
\kring \fs^{\ell/ (\ell -1)}.
$ 
Recall that $\triv$ is the constant $\kring \fb\op$-module $\kring$, which identifies as $\bigoplus_{\ell \in \nat} \triv_\ell$. 
Using this, one has the isomorphism
\[
\primfs \odot_{\fb\op} \triv \cong \bigoplus_{\ell \in \nat} \primfs \odot_{\fb\op} \triv_\ell.
\]

Theorem \ref{thm:coker_theta_ell} thus has the following immediate consequence, showing how $\kring \fs$ is built from $\primfs$ as a $\kring \fb$-bimodule (aka. $\kring \fb \otimes \kring \fb\op$-module):

\begin{corollary}
\label{cor:kring_fs}
There is an isomorphism of $\kring \fb \otimes \kring \fb\op$-modules
\[
\primfs \odot_{\fb\op} \triv 
\cong 
\kring \fs 
\oplus 
\bigoplus _{ \substack{a \in \nat \\ 0< \ell \in \nat}}  \big(\sgn_a \boxtimes S_{(\ell, 1^a)} \big).
\]
\end{corollary}

\section{Calculating $\primfs$ as a $\kring \fb$-bimodule}
\label{sect:primfs}

In this section, we  use Corollary \ref{cor:kring_fs} to calculate the underlying $\kring \fb$-bimodule of $\primfs$ in terms of $\kring \fs$ (see  Theorem \ref{thm:primfs_bimodule}).  This then feeds back into Theorem \ref{thm:coker_theta_ell} so as to yield Corollary \ref{cor:subquotients}, which determines the underlying $\kring \fb$-bimodule of all the subquotients of the filtration of $\kring \fs$.

\begin{remark}
\label{rem:compare_Theorem3_PV2}
Theorem \ref{thm:primfs_bimodule}  recovers  \cite[Theorem 3]{MR4518761}, which was established by using  different methods. 
 Some explanation is necessary: by  \cite[Corollary 5.2 and Theorem 5.4]{MR4518761}, the result \cite[Theorem 3]{MR4518761} corresponds to the calculation of $H^0 (\mathbb{C}(\mathbf{b}, \mathbf{a}))$, which identifies with $\primfs (\mathbf{b}, \mathbf{a})$, as observed in Remark \ref{rem:primfs_explicit}. 

The approach in  \cite{MR4518761} relies on calculating the cohomology of the complex $\mathbb{C}(\mathbf{b}, \mathbf{a})$; this complex is constructed using the $\kring \finj\op$-module structure of $\kring \fs (-, \mathbf{a})$. The positive degree cohomology of $\mathbb{C}(\mathbf{b}, \mathbf{a})$ is sufficiently sparse and calculable that an Euler-Poincaré characteristic argument allows $H^0 (\mathbb{C}(\mathbf{b}, \mathbf{a}))$ to be calculated. Most of this argument can be carried out working over an arbitrary commutative ring. 
\end{remark}

\subsection{Recovering information from   $M \odot_{\fb\op} \triv$}

The proof of Theorem \ref{thm:primfs_bimodule} is based on the fact that one can recover the isomorphism class of a $\kring \fb\op$-module $M$ from $M \odot_{\fb\op} \triv$, as stated explicitly in Proposition \ref{prop:invert}. (This technique is also used in \cite{P_finmod} (see \cite[Section 8]{P_finmod}); the details are outlined below so as to keep the exposition self-contained.)

We  work with isomorphism classes of objects and with virtual objects; to do so, we pass to the Grothendieck group of $\kring \fb\op$-modules. 

\begin{remark}
We only work with $\kring \fb\op$-modules that take finite-dimensional values. Thus, after evaluation on some finite set, all calculations reduce to studying  finite objects. For simplicity of exposition, such details have been suppressed. 
\end{remark}

\begin{notation}
\label{nota:groth_gp}
\ 
\begin{enumerate}
\item 
For a $\kring \fb\op$-module $M$, its class in the Grothendieck group is denoted $[M]$. 
\item 
For $\kring \fb\op$-modules $M$, $N$, write $[M] \odot_{\fb\op} [N]$ for the element of the Grothendieck group $[M \odot_{\fb\op} N]$; this only depends on the classes $[M]$, $[N]$.
\end{enumerate}
\end{notation}

The following observation ensures that the infinite sum appearing in Proposition \ref{prop:invert} only has finitely-many non-zero terms when evaluated on any $\mathbf{n}$, hence is well-defined.  

\begin{lemma}
\label{lem:vanishing}
For $N$ a $\kring \fb\op$-module and $t\in \nat$, $(N \odot_{\fb\op} \sgn_t) (\mathbf{n}) =0$ for $t>n$.
\end{lemma}

The following result shows how to recover $[M]$ from $[M \odot_{\fb\op} \triv]$:

\begin{proposition}
\label{prop:invert}
For $M$ a $\kring \fb\op$-module, one has the equality:
\[
[M] = \sum_{t \in \nat} (-1)^t [M \odot_{\fb\op} \triv] \odot_{\fb\op}[\sgn_t]
= 
[M \odot_{\fb\op}\triv] \odot_{\fb\op}\Big( 
\sum_{t \in \nat} (-1)^t  [\sgn_t]
\Big) 
.
\]
\end{proposition}

\begin{proof}
By associativity of $\odot_{\fb\op}$, one reduces to the case $M = \kring$ supported on $\mathbf{0}$.

The equality
$ 
[\kring] = \sum_{t \in \nat} (-1)^t [\triv] \odot_{\fb\op}[\sgn_t]
$  
is seen as follows. Evaluated on $\mathbf{0}$ one recovers $[\kring]$ from the right hand side. Thus, to conclude, for any $n>0$, we require to show that 
\begin{eqnarray}
\label{eqn:de_Rham}
\sum_{t=0}^n (-1)^t [\triv_{n-t}] \odot_{\fb\op}[\sgn_t] =0.
\end{eqnarray}
This can be checked by using the Pieri formula to calculate the $\kring \sym_n\op$-representations $\triv_{n-t} \odot_{\fb\op}\sgn_t$ (cf. the proof of Lemma \ref{lem:identify_Lambda^a_pbar}).
\end{proof}

\begin{remark}
Proposition \ref{prop:invert} can be paraphrased as stating that the functor $-\odot_{\fb\op} \triv$ is invertible.
\end{remark}

\subsection{Recovering $\primfs$}

Proposition \ref{prop:invert} applies using Corollary \ref{cor:kring_fs} to calculate the isomorphism class of the underlying $\kring \fb$-bimodule of $\primfs$:

\begin{theorem}
\label{thm:primfs_bimodule}
In the Grothendieck group of $\kring \fb \otimes \kring \fb\op$-modules, one has the equality
\begin{eqnarray}
\label{eqn:primfs}
[\primfs] + \sum_{b> a \geq 0} (-1)^{b-a} [\sgn_a \boxtimes \sgn_b] 
= 
[\kring \fs] \odot_{\fb\op}\Big( 
\sum_{t \in \nat} (-1)^t  [\sgn_t]
\Big) .
\end{eqnarray}
\end{theorem}

\begin{proof}
By Corollary \ref{cor:kring_fs}, one has the equality in the Grothendieck group of $\kring \fb$-bimodules:
\[
[\primfs \odot_{\fb\op}\triv] 
=
[\kring \fs] 
+ 
\sum_{ \substack{a \in \nat \\ 0< \ell \in \nat}} [\sgn_a \boxtimes S_{(\ell, 1^a)}].
\]
Since $\odot_{\fb\op}$ is symmetric, 
Proposition \ref{prop:invert} gives the equality 
\[
[\kring \fs] \odot_{\fb\op}\Big( 
\sum_{t \in \nat} (-1)^t  [\sgn_t]
\Big)
\odot_{\fb\op}[\triv] 
= 
[\kring \fs].
\]

To conclude, we require to show that 
\[
\sum_{ \substack{a \in \nat \\ 0< \ell \in \nat}} [\sgn_a \boxtimes S_{(\ell, 1^a)}]
+ 
 \sum_{b> a \geq 0} (-1)^{b-a} [\sgn_a \boxtimes \sgn_b] \odot_{\fb\op}[\triv] =0.
\]
Equivalently, for each  $\ell >0$ and $a \in \nat$, we require to show that  
\[
\mathcal{V}(\ell, a):= 
\ 
[S_{(\ell, 1^a)}] 
+ 
\sum_{b>a}
(-1)^{b-a}
[\sgn_b \odot_{\fb\op}\triv_{\ell + a -b}] ,
\]
is zero,  since we can neglect the common factor of $\sgn_a \boxtimes -$.

This is proved by increasing induction upon $a \geq 0$, starting from the case $a=0$, for which $\mathcal{V} (\ell, 0) =0$ by 
(\ref{eqn:de_Rham}). One has the identity $\sgn_{a+1} \odot_{\fb\op} \triv_{\ell -1} \cong S_{(\ell, 1^a)} \oplus S_{(\ell -1, 1^{a+1})}$ (where, for $\ell=1$, the second term is understood to be zero) by Pieri's formula (cf. the proof of Lemma \ref{lem:identify_Lambda^a_pbar}), which gives the equality $\mathcal{V} (\ell, a) = - \mathcal{V} (\ell, a+1)$. It follows inductively that $\mathcal{V}(\ell, a)$ vanishes for all $a$. 
\end{proof}

\subsection{Identifying the higher subquotients of the filtration of $\kring \fs$}
Theorem \ref{thm:coker_theta_ell} yields the equality in the Grothendieck group of $\kring \fb \otimes \kring \fb\op$-modules  
for $\ell >0$
\begin{eqnarray}
\label{eqn:kkfs_ell/ell-1}
\ 
[\kring \fs^{\ell/(\ell -1)}] 
=
[\primfs] \odot_{\fb\op}[\triv_\ell]
- 
\sum_{a \in \nat}
[\sgn_a \boxtimes S_{(\ell, 1^a)}].
\end{eqnarray}
Hence, using Theorem \ref{thm:primfs_bimodule}, one can identify $[\kring \fs^{\ell/(\ell -1)}] $ in terms of $\kring \fs$:

\begin{corollary}
\label{cor:subquotients}
For $0< \ell \in \nat$, in the Grothendieck group of $\kring\fb \otimes \kring \fb\op$-modules, there is an equality
\begin{eqnarray*}
[\kring \fs^{\ell/(\ell -1)}] 
&=&
[\kring \fs] \odot\Big( 
\sum_{t \in \nat} (-1)^t  [\triv_\ell \odot\sgn_t]
\Big)
- 
\sum_{b>a\geq 0} (-1)^{b-a} [\sgn_a \boxtimes (\sgn_b \odot\triv_\ell)] 
\\
&& 
- 
\sum _{a \in \nat} [\sgn_a \boxtimes S_{(\ell, 1^a)}],
\end{eqnarray*}
in which $\odot$ denotes $\odot_{\fb\op}$.
\end{corollary}

\subsection*{Acknowledgement}

This project was inspired in part by joint work with Christine Vespa; the author is grateful for her ongoing interest. The author also thanks Vladimir Dotsenko for his invitation to Strasbourg in January 2024 (financed by his {\em IUF}) and for discussions related to this project. The author is also very grateful to two anonymous referees for their  constructive comments.


\EditInfo{July 17, 2025}{November 13, 2025}{Roozbeh Hazrat}


{\small
    \begin{thebibliography}{10}
    
        \bibitem{MR2350124}
        B.~Feigin and B.~Shoikhet.
        \newblock On {$[A,A]/[A,[A,A]]$} and on a {$W_n$}-action on the consecutive
        commutators of free associative algebras.
        \newblock {\em Math. Res. Lett.}, 14(5):781--795, 2007.
        
        \bibitem{MR1153249}
        W.~Fulton and J.~Harris.
        \newblock {\em Representation theory}, volume 129 of {\em Graduate Texts in
        Mathematics}.
        \newblock Springer-Verlag, New York, 1991.
        \newblock A first course, Readings in Mathematics.
        
        \bibitem{Phh}
        T.~Pirashvili.
        \newblock Hodge decomposition for higher order {H}ochschild homology.
        \newblock {\em Ann. Sci. \'Ecole Norm. Sup. (4)}, 33(2):151--179, 2000.
        
        \bibitem{MR4398644}
        L.~Positselski.
        \newblock {\em Relative nonhomogeneous {K}oszul duality}.
        \newblock Frontiers in Mathematics. Birkh\"{a}user/Springer, Cham, [2021]
        \copyright 2021.
        
        \bibitem{P_finmod}
        G.~{Powell}.
        \newblock {Functors on the category of finite sets revisited}.
        \newblock {\em arXiv:2407.11623}, July 2024.
        \newblock Accepted, {\em Ann. Math. Blaise Pascal}.
        
        \bibitem{P_rel_kos}
        G.~{Powell}.
        \newblock {Relative nonhomogeneous Koszul duality for PROPs associated to
        nonaugmented operads}.
        \newblock {\em arXiv:2406.08132}, June 2024.
        
        \bibitem{MR4518761}
        G.~Powell and C.~Vespa.
        \newblock A {P}irashvili-type theorem for functors on non-empty finite sets.
        \newblock {\em Glasg. Math. J.}, 65(1):1--61, 2023.
        
        \bibitem{PV}
        G.~Powell and C.~Vespa.
        \newblock Higher {H}ochschild homology and exponential functors.
        \newblock {\em Bull. Soc. Math. France}, 153(1):1--141, 2025.
        
        \bibitem{MR3430359}
        S.~V. Sam and A.~Snowden.
        \newblock G{L}-equivariant modules over polynomial rings in infinitely many
        variables.
        \newblock {\em Trans. Amer. Math. Soc.}, 368(2):1097--1158, 2016.
        
        \bibitem{MR3556290}
        S.~V. Sam and A.~Snowden.
        \newblock Gr\"{o}bner methods for representations of combinatorial categories.
        \newblock {\em J. Amer. Math. Soc.}, 30(1):159--203, 2017.
        
        \bibitem{MR4788713}
        S.~V. Sam, A.~Snowden, and P.~Tosteson.
        \newblock Polynomial representations of the {W}itt {L}ie algebra.
        \newblock {\em Int. Math. Res. Not. IMRN}, (16):11688--11710, 2024.
        
        \bibitem{MR4546495}
        P.~Tosteson.
        \newblock Categorifications of rational {H}ilbert series and characters of
        {$FS^{\rm op}$} modules.
        \newblock {\em Algebra Number Theory}, 16(10):2433--2491, 2022.
        
        \bibitem{MR3982870}
        V.~Turchin and T.~Willwacher.
        \newblock Hochschild-{P}irashvili homology on suspensions and representations
        of {${\rm Out}(F_n)$}.
        \newblock {\em Ann. Sci. \'{E}c. Norm. Sup\'{e}r. (4)}, 52(3):761--795, 2019.
        
        \bibitem{2014arXiv1406.0786W}
        J.~D. {Wiltshire-Gordon}.
        \newblock {Uniformly Presented Vector Spaces}.
        \newblock {\em arXiv:1406.0786}, June 2014.

    \end{thebibliography}
}
\end{document}